# Is the study of Indigenous mathematics ill-directed or beneficial?

Hongzhang Xu[1], Mathematical Sciences Institute & Fenner School of Environment and Society,
and
Rowena Ball[2], Mathematical Sciences Institute,

Australian National University, Canberra ACT 2601 Australia

**Abstract**

The 'old lie' of mathematical inadequacy of Indigenous communities has been curiously persistent despite increasing evidence shows that many Indigenous communities practiced mathematics. Attempts to study and teach Indigenous mathematical knowledge have always been questioned and even denied validity. The Aboriginal and Torres Strait Islander Histories and Cultures cross-curriculum priority in the F-10 Australian schools curriculum, from 2022 onwards, includes content elaborations related to Indigenous mathematics, which have been developed and refined by expert Indigenous advisers. We celebrate this initiative, but experience also tells us to expect some resistance from sectors of the education communities who hold to an exclusively Anglo-European provenance of mathematics. Through this review article we seek to constructively forestall potential pushback and address concerns regarding the legitimacy and pedagogical value of Indigenous mathematics, by countering with evidence some published claims of mathematical inadequacies of Australian First Nations cultures.

## 1. Introduction

Indigenous sciences, in all their diversity, depth and vitality, are actively developed by First Nations practitioners and increasingly respected by mainstream scientists. Global recognition of and advocacy for Indigenous sciences and technology have enriched fields such as knowledge of ecology (Bardsley et al., 2019), geography (Ioris et al., 2019), geology (Larkin et al., 2012), astronomy (Hamacher et al., 2020), hydrology (Moggridge et al., 2022) and climate change (Leal Filho et al., 2022; Whyte, 2018). As scientific approaches to nature, technologies and human society seem cross-culturally universal, thinking and communicating mathematically also is likely to be a universally human practice (Verran, 2000).

Yet the study and teaching of Indigenous mathematical knowledge have always been called into question and even denied validity. Ethnomathematics has been denigrated in international opinion media as a 'fringe absurdity' (Brooks, 2021). In response to a (then) draft proposal to include Indigenous mathematical perspectives in a revised Australian national school curriculum Deakin (2010) argued that 'the topic of Indigenous mathematics does not exist…attempts to discover an Indigenous Mathematics are…**ill-directed**…There is no indigenous tradition of mathematics, properly so-called, in this country [Australia]'.

[1] Corresponding author, email: hongzhang.xu@anu.edu.au, https://orcid.org/0000-0001-8904-2976
[2] Email: rowena.ball@anu.edu.au, http://orcid.org/0000-0002-3551-3012

Michael Deakin (1939-2014) was an Australian mathematician and mathematics educator. Although more than a decade has passed since publication of this article (Deakin, 2010), the views expressed therein have been reiterated for generations (Thomas, 1996) and have persisted in mainstream mathematics, by default or indifference if not actively, and as a result have been reflected in school and university mathematics curricula, activities and positions taken by mathematics professional bodies, and educational policy and ideology.

In this article we advance the epistemology of both Indigenous and non-Indigenous mathematics through our replies to Deakin's five main critiques, which are tabulated in section 2. In doing so, we review some relevant literature on Indigenous mathematics. We refer to the critiques generically as 'the Gatekeeper critiques' for three reasons:

1) They are representative expressions which reflect a dogma of exclusiveness of mathematics to an Anglo-European tradition that pervades the community to this day;

2) This dogma is defended actively by guardians, or gatekeepers, against any perceived illegitimate incursions or subversions;

3) The 'gatekeeping method' of controlling access and inclusion is used to disallow the legitimacy of Indigenous mathematical knowledge.

In section 3 we examine the Gatekeeper critiques critically and in section 4 we propose a new rapport whereby Indigenous and non-Indigenous mathematics may inform each other and advance together.

## 2. The Gatekeeper critiques of Australian First Nations mathematics

In pushback against evidence that Australian First Nations peoples practised mathematics, Deakin (2010) questions and dismisses or impugns several reports on Indigenous mathematical knowledge, concluding (unsurprisingly) that the examples described most definitely are *not* mathematics.

What *is* surprising – or perhaps not – is that the Gatekeeper critiques have not to date been challenged in the literature. Among our mathematical colleagues, many if not most recognize that there are well-documented non-Eurocentric traditions of mathematics, such as Mayan, Chinese and Indian mathematics (Joseph, 2009; Peat, 2002; Tiles, 2002). For some, especially global gatekeepers of pure mathematics, however, the notion of 'Indigenous mathematics' as systems of knowledge utilising symbolism, abstraction, classification, and quantitation is going a step too far.

Many other international mainstream mathematicians also have questioned the existence of Indigenous mathematics (Harris, 1987; Millroy, 1991; Pais, 2011). No written language, no coherent system of numerals and rough measurement are the main reasons given. Yves Chevallard (1990) denied Palus Gerdes' finding that mathematics comes from inside of African, Asian and American-Indian cultures. Chevallard (1990) called indigenous physics and mathematics as proto-physics and proto-mathematics, respectively, because of the inadequacy of 'fully-fledged sciences' in indigenous communities (Chevallard, 1990). Rowlands and Carson (2002) regards mathematics of the people as non-rigorous and informal, being distinctly different from mathematics of 'elite', i.e., academic or pure mathematics. Pais (2011) proposes that school should only provide students with decontextualized universal knowledge by excluding cultural components (Pais, 2011).

The Gatekeeper critiques of Australian First Nations mathematics collecting most queries on indigenous mathematics across the world. We present the Gatekeeper critiques in Table 1.

Table 1 The Gatekeeper critiques (Deakin, 2010). (Page numbers in the table refer to those of the cited article.)

| Critiques | Quotes |
| --- | --- |
| 1. Australian First Nations had no scripted languages and therefore no mathematics. | p.234, quoted from Cohen (2006): "*There were no written languages…and probably for that reason there was no systematic development of mathematics*" |
| 2. Australian First Nations' numeral systems lacked capability to carry out arithmetic. Systems of numerals are not mathematics. | p.234, quoted from Rudder (1999): "*They do classification almost completely without numbers…[numerical] concepts are not precise enough to lend themselves even to the simplest arithmetic manipulation*". p.235: Published evidence on Indigenous number systems is dismissed by gainsay: "*the existence of a coherent system of numerals can hardly be said to constitute mathematics*". |
| 3. Knowledge systems and practices such seasonal calendars, time of day and astronomy involve no mathematics and have no relation to measurement and geometry. | p.235: "*calendric time and time of day…is mathematical in only the very broadest sense of that word*" p.236: "*counting seasons is not mathematics, and quite how it relates to measurement and geometry is anybody's guess*" p.236: "*Aboriginal recognition [of] solar eclipses…is really glaringly obvious and furthermore involves no mathematics…Were there an Aboriginal method of predicting eclipses the case [for the existence of Aboriginal mathematics] would be made*" |
| 4. There is no mathematics in Aboriginal and Torres Strait Islander groups' kinship rules | Dismissal by *non sequitur*: Western mathematics does not describe Indigenous kinship rules perfectly, therefore they cannot be mathematical. p.235: "*group-theoretical discourse…is not a product of Aboriginal culture…group theory comes from mainstream western tradition*". |
| 5. Imprecise, inadequate and inconsistent Indigenous mathematical systems were inferior and now are irrelevant. | p.236: "*Our mainstream culture…needs the precision of modern science. Earlier cultures…are nowadays superseded by [western science]*. |

## 3. Reply to the Gatekeeper critiques

Some denigratory statements on Australian Indigenous mathematics have been successfully rebated by John Harris (1987) and Judith Stokes (1982). However, the existence of mathematics hold by the Australian First Nations people has still been rejected by shifting the focus to other spheres that have been little studied. There is always more to be questioned and the Gatekeeper critiques may never be stopped. A powerful rebuttal is needed to challenge their beliefs on the origin and superiority of mainstream mathematics. To do so, we broaden these critiques to refute past, existing and future deniers in pushback against evidence that First Nations peoples practised mathematics in Australia and the rest of the world. We include diverse global indigenous cases to reply to these representative critiques on mathematical inadequacy of First Nations people.

Deakin's argument is deductive and follows a top-down approach in which his conclusion is based on a premise that is assumed to be true (Weddle, 1979). The premise relies on answers to the questions: What is mathematics? What is included in mathematics? Deakin does not answer these questions but assumes that the reader knows and agrees with him. This lack of clear premise makes his critiques dubious and subjective at the outset. In this section, we reply to the Gatekeeper critiques in Table 1 by reviewing some history of mathematics and latest research progress in Indigenous mathematics. We refer to documented mathematical knowledge from Indigenous and historical cultures worldwide.

*Critique 1. No written language, ergo, no mathematics*

It is often assumed that mathematics cannot be done without a system of written language expression. This is because mathematics is believed to be an abstract, context-free pastime tied to a formal system which relies upon canonical sets of symbols that have to be written down. This idea has led to a view of mathematics are existing in the æther, outside ordinary human activity and devoid of social and cultural considerations (Millroy, 1991).

However, it is often forgotten or ignored that mathematical expertise, for example in arithmetic and geometry, is originally developed to solve real problems that arise in daily life (Ascher & Ascher, 1986; Sizer, 2000). In any culture there is an agreed, structured and self-consistent system to classify, transform, quantify, communicate and predict patterns and cycles of importance, whether it be in unwritten or written forms.

Much of ordinary day-to-day arithmetic and geometry performed by 'illiterate' women, artisans, carpenters and many other workers are unwritten and even unspoken (Wood, 2000). The apprentice learns by watching carefully then doing the mathematics themselves. The use of tools–an unwritten approach–to support arithmetic has a long history; there are different media for recording and computing with numbers, including stones, twigs, knots and notches (Hansson, 2018). People of many Indigenous Pacific and Australian nations can use parts of the body to count quickly and accurately (Goetzfridt, 2007; Owens & Lean, 2018; Wood, 2000), communicating methods, operations and results through speaking, listening and gesture. Weaving skills were taught unwritten to next generations to construct the numerical relationships that give rise to the desired complex geometrical designs with symmetries (Hansson, 2018). Knotted *quipus* were used by 'illiterate' Inca people of South American Andes regions to allot land and levy taxes (Ascher & Ascher, 2013). The quipu (Figure 1), with its columns of base-10 numerical data encoded as knots, can be thought of as a spreadsheet, and it seems likely that the Inca knew and applied some array and matrix operations.

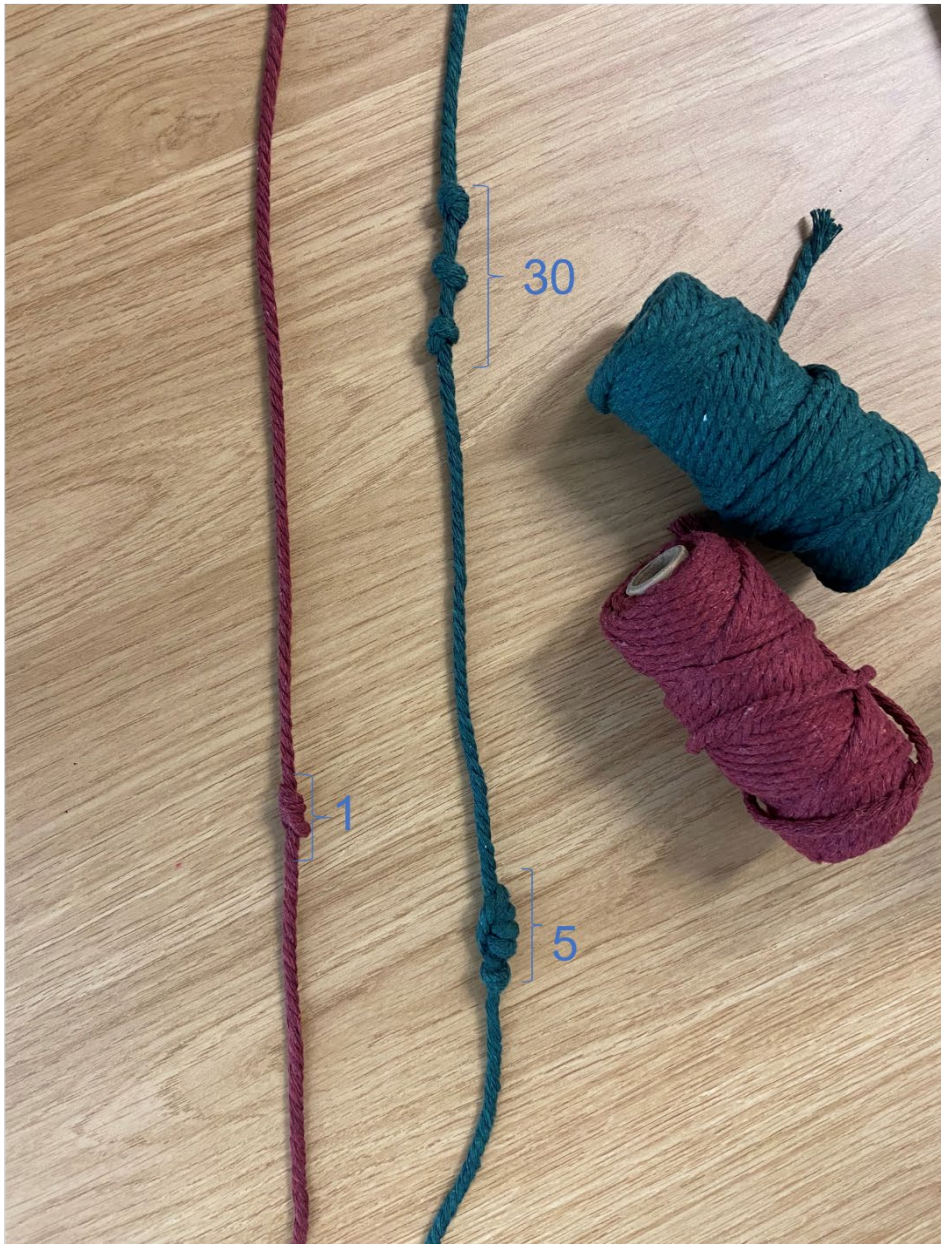

Figure 1. Two basic units of a quipu showing numbers 1 and 35. (A work-in-progress by Rowena Ball). © Authors

*Dan*, an Indigenous language of central Liberia, is non-written but *Dan* speakers can carry out arithmetic operations orally, including addition, subtraction and division, play games that require fast counting, tracking and calculating skills, and practice geometric principles in constructing buildings (Sternstein, 2008). Fractal geometry, developed to a high art in Western mathematics from the late 1960s and executed *in silico*, has non-Western antecedents that were implemented in the built environment in Africa (Eglash, 1998). Chaology and fractal geometry have also been a part of traditional Chinese architectural and garden design for thousands of years (Li & Liao, 1998).

Moreover, Latin-European mathematics – which became mainstream – originated and was carried out in unwritten forms. Verbal or non-written approaches played important roles in ancient Greek mathematics, which relied on diagrams, sculpture and verbal echoes. References and expansions were chanted orally (Netz, 2005).

We need to understand that mathematics is socially constructed in the context of a community, where meaning is negotiated, and conventions are agreed upon certain groups of people. Writing them down is the dominant approach today from Western mathematics but it is not the only one.

*Critique 2. No coherent system of numerals*

Pursuant to the above discussion, we have diverse and well-documented evidence that societies nominally without writing systems have well-developed and useful numeral systems. Is that mathematics? Not so, according to Deakin (2010, p. 234), who uses the terms 'proto-numerals' and 'proto-mathematics' to claim that Australian First Nations numeral systems are not mathematics because they cannot support basic arithmetic manipulation.

In fact the 'turtle egg mathematics' practised by Yolŋu people of north east Arnhem land, selected by Deakin as an example of imprecise or ambiguous use of number counting, gives perfectly accurate results; rather, it appears that Deakin's understanding of this system is inaccurate. Cooke (1991, p. 12) and Lloyd et al. (2016, p. 7) explain the word '*rulu*', which means a group or bundle in general. However, prefix and suffix will be added to scope down its meaning to a specific number. For example, '*wangany-rulu*' means five. A closer look at the Yolŋu turtle eggs sharing strategies finds that subtractive division in base 5 is the method applied, as illustrated in Figure 2.

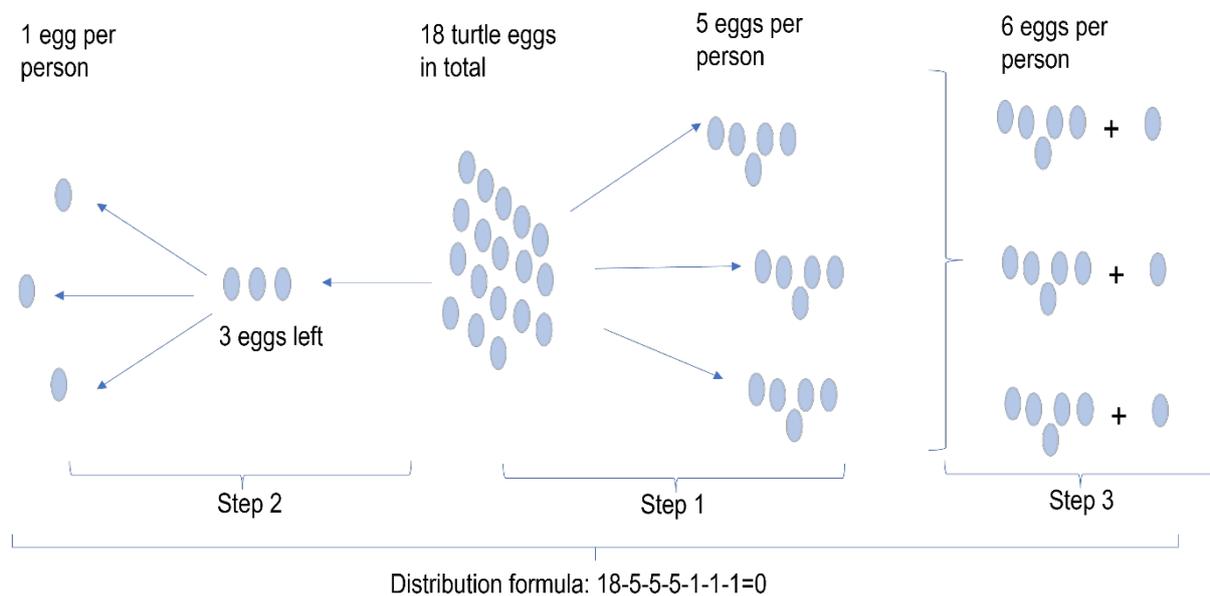

Figure 2. Visualization of Yolŋu People's turtle eggs distribution process in a base 5 numeral system

We cite also the following examples that are well-documented:

- The Anindilyakwa people of Groote Eylandt, Northern Territory, use systematic numerals up to 20 (Harris, 1987), i.e., a base 20 system.

- Glendon Lean's research on nearly 900 counting systems of Papua New Guinea, Oceania and Irian Jaya provides a wide range of counting and tallying systems with diverse base systems, from base 2 to base 27 (Owens, 2001).
- Research on body-part tally systems in Papua New Guinea finds more counting systems, ranging from base 27 to 35 (Dwyer & Minnegal, 2016) and in base 10 with base 2 as reflection (Owens & Muke, 2020).

In addition, much of mathematics is non-numeric. Numbers are not necessary in classification, which refers to the arrangement or sorting of objects into groups based on their common property. Many existing and widely used classification systems are qualitative or semi-qualitative. For example, high-medium-low, A-B-C-D and Yes-No are not numeric. Non-western systematic classification systems were developed thousands of years ago by Buddhism (e.g., six paths) and School of Naturalists (e.g., wuxing and bagua) (Gethin, 1998; Law & Kesti, 2014). Indigenous communities in Australia have detailed and accurate taxonomy systems for plants and animals (Neale & Kelly, 2020).

*Critique 3. No measurement of time and space*

Although Deakin (2010, p. 236) acknowledges that 'astronomy has historically been seen as a branch of mathematics', he shuts the gate against Indigenous knowledge on astronomy by claiming it makes no measurements of time and space and therefore is not mathematical. Such misconceptions are still prevalent, and have been challenged by (Hamacher & Norris, 2011).

The connections between Indigenous astronomy and mathematics are becoming clearer. Much richer Aboriginal and Torres Strait Islander knowledge of astronomy has been revealed, *through* what Deakin (2010, p. 236) refers to as the 'elaboration of stories' (Noon & Napoli, 2022). Indigenous knowledge on eclipses has shown an insightful understanding of the earth-moon-sun system and its correlations with terrestrial events, including seasonal cycles and weather forecasting (Bureau of Meteorology, 2022; Hamacher et al., 2019; Norris & Hamacher, 2009). For example, the Yolŋu people of Arnhem Land always have understood the relationship between the Moon cycles and the ocean tides. Yolŋu Elders can predict the time and height of the next tide after seeing the position and phase of the Moon, which contrasts with Galileo's incorrect explanation of the tides (Norris, 2016).

But Deakin (2010, p. 236) has devised an infallible test for the existence of Indigenous mathematics! This is that there must be 'an Aboriginal method of predicting eclipse'. To predict an eclipse, one needs clear and accurate understanding of the relationships between the motions of the Sun and Moon. In spite of the challenge, the answer is yes. Hamacher & Norris (2011) report a prediction by Aboriginal people of a solar eclipse that occurred on 22 November 1900, which was described in a letter dated in December 1899.

Spatial coordinate systems are well-defined in First Nations societies. The use of cardinal directions (north, south, east and west) is common amongst Aboriginal language groups in Australia (Norris & Harney, 2014). Guugu Yimithirr language speakers use cardinal directions rather than concepts of left, right, behind or front to determine objects spatially and in time (Norris, 2016). Thus the mathematical concept of vectors is strongly embedded culturally. Strong directionality and astronomical knowledge enabled Aboriginal people to undertake long-distance travel and develop extensive trading networks (Forster, 2021; Norris & Harney, 2014).

It is questionable to simply deny the existence of mathematics by ignoring, or attributing to chance, the mathematical thought and practice evidenced by these capacities.

*Critique 4. No mathematics in kinship rules*

Research on Aboriginal kinship rules is one of first bridges between Indigenous and mainstream mathematics. The mathematics of Aboriginal kinship rules was first described by Mathews (1900), and many studies since have found the existence of group theoretical knowledge in such systems (Keen, 1988; Laush, 1980; Radcliffe-Brown, 1930). Broadly, group theory is a suite of largely non-numerical mathematics involving canonical classifications of sets of objects with mathematical properties. A probabilistic approach by Field (2021) has revealed the intergenerational social trade-offs and benefits of the Gamilaraay kinship system, and the simpler aspects of the mathematics of Yolŋu kinship has been elucidated for a teaching audience by Matthews (2020).

However, Deakin denies connections between group theory and kinship rules by stating that group theory was invented by Western mathematicians, and selectively claims that mathematical studies of kinship rules are just an application of Western theory (Figure 3). These algebraic analyses of kinship Cargal (1978) are bridges between Aboriginal knowledge system and modern mathematics but are misinterpreted by Deakin (p. 235) as 'pseudo-mathematics'.

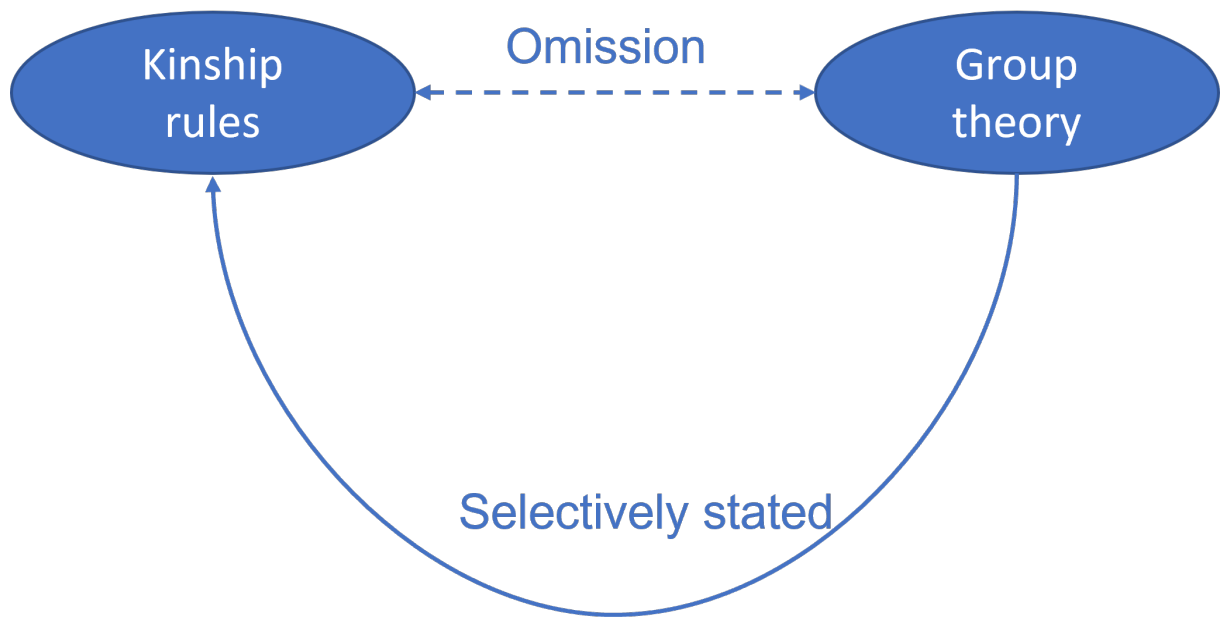

Figure 3. Deakin's misleading argument

The group-theoretical knowledge manifested by Aboriginal kinship rules has been applied for thousands of years. In 1949, André Weil used group-theoretical techniques to model Australian marriage systems, marking the beginning of mathematical anthropology (Ballonoff, 2018; Rauff, 2016).

Not all kinship rules have a group theoretical interpretation though, as indicated in Figure 3. Such irregularities are the 'doubts' and 'exceptions' mentioned by Deakin. But Aboriginal knowledge is an integrated system not a discipline-segregated one such as mainstream mathematics (Sizer, 2000). Marriage and other kinship activities need mathematics but cannot be done by using only mathematics. Mathematics is likely to benefit from studying the irregularities to understand their possible relations (Figure 4), rather than gatekeeping to deny identified connections. Kinship rules are a clear indication that First Nations societies mastered the powerful instruments of mathematics to live and manage the community better.

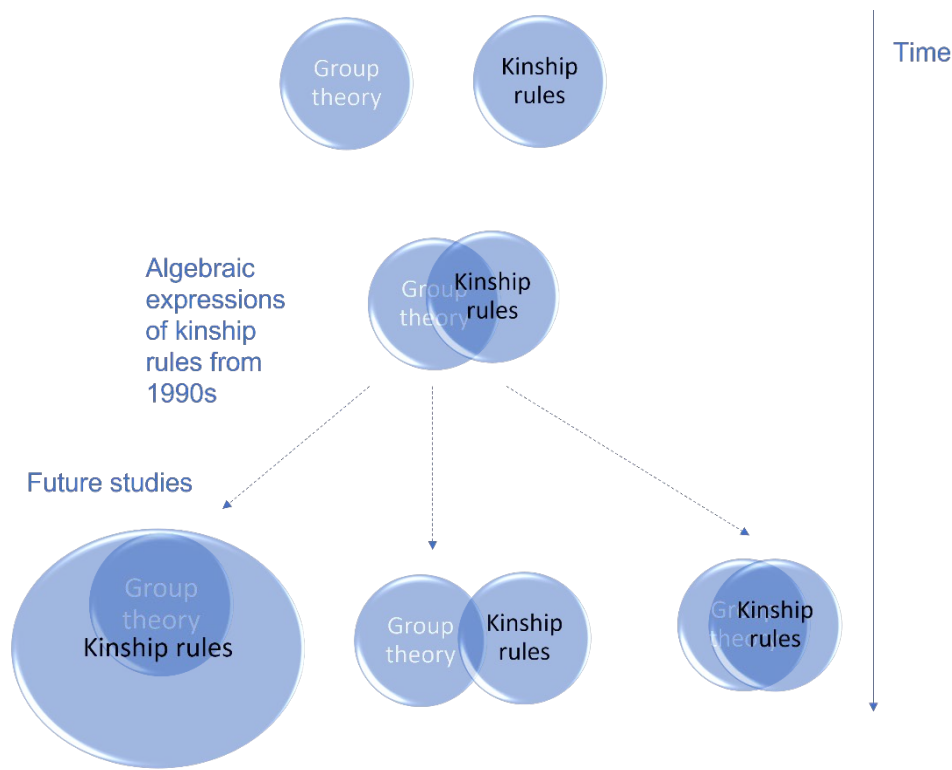

Figure 4. Research progress of group theory and kinships rules

*Critique 5. Inability of conducting complex weather and season change prediction*

Gatekeeping is again the main forte of the fifth critique, which relies on the belief that only Western culture has developed an advanced system to understand weather patterns and season change. However, great uncertainties still exist in numerical long-range weather and climate prediction because of the existence of 'chaos' in dynamic systems (Aurell et al., 1997). It is thus that Gatekeepers assert that Indigenous People are not able to conduct complex prediction due to their outdated knowledge system.

First Nations people developed unique methods of living which enabled life and community to flourish in even the harshest environments and have allowed the sustained development of the world's oldest continuing culture. Seasonal knowledge was important for survival – many First Nations peoples changed location seasonally to optimize water availability, and for associated hunting, fishing, plant-food gathering and cultivation (Hamacher et al., 2019). Being examined for tens of thousands of years, Indigenous meteorology is a precious and irreplaceable heritage. It is an intimate knowledge of plant and animal cycles and contains details of the intricate connections in the natural world (Bureau of Meteorology, 2022).

Indigenous meteorology is a precious and irreplaceable heritage. It relies on observational knowledge of the skies and plant and animal cycles (Green et al., 2010). For example, stellar scintillation has been used by Torres Strait Islanders to predict weather and seasonal change (Hamacher et al., 2019). Scintillating stars that appear sharp to the eye indicate dry, clear skies. Humidity and ice crystals in the atmosphere can blur the stars and warn of a coming storm or cyclone. Thus, the jagged appearance of scintillating stars denotes moisture, thus warning of a coming storm or cyclone (Hamacher et al., 2019). When most of the visible stars are blue, it indicates moisture, as the red and green wavelengths of light are absorbed by the humidity. In the Torres Strait, twinkling blue stars indicate hot weather with increasing humidity as the wet monsoon season approaches. The blurry characteristic of stars during periods of high humidity is known to many farmers (Hamacher et al., 2019). In addition, as we introduced above, Yolŋu People can predict tides accurately according to their knowledge of the environment and information gathering from the stellar system (Norris, 2016). Aboriginal people can also predict weather and season change based on biological indicators. For example, Torres Strait Islanders can predict windy and rainy days by observing the movement of sharks (*Beizam*) (Green et al., 2010). Walabunnba people knows that there will be a lot of rain when hearing the mirrlarr (rain bird) calls out (Bureau of Meteorology, 2022). The flowering of the boo'kerrikin (Acacia decurrens) is an indication for the D'harawal people to know the start of gentle spring rains and end of the cold and windy weather (Bureau of Meteorology, 2022).

Indigenous understandings, measurement and predictions of the weather and season change are complicate processes (Clarke, 2009). Different streams are included, compared and synthesised to predict the change of weather and seasons. Rather than denying and ignoring the merits of First Nations meteorology and science, co-production of weather and season knowledge with Aboriginal people could complement and may add new insights into our existing forecasting system (Bureau of Meteorology, 2022). Will the algorithm be improved, and computing speeded be accelerated by embedding First Nations knowledge? May modern computational weather prediction be improved by embedding First Nations knowledge? Could non-meteorological observations, such as plants and animals' behavior change, been used to improve existing weather models? These potentials cannot be realized if we disrespect the knowledges of First Nations cultures.

## 4. Discussion and conclusion

Our responses above indicate that the Gatekeeper critiques on Indigenous mathematical inadequacy are not rooted in either sound evidence or clear logic. We also illustrate that mathematics produced by Indigenous People can contribute to the economic and technological development of our current 'modern' world.

The mathematical ideas of Indigenous people have been disregarded as part of the colonial strategy of gatekeeping the periphery of structured scientific knowledge (Ball, 2015; Goetzfridt, 2007). For example, medieval Europe justified the 'success' of their conquest of the Americas due to their cultural superiority despite Mayans were more developed in many fields of mathematics and sciences than the Europe in 1500s (Rosa & Orey, 2008). To gatekeep the cultural superiority, European invaders destroyed the Mayan libraries (Rosa & Orey, 2008). But the surviving texts left reveal the sophistication of Mayan knowledge of astronomy and mathematics, especially their advanced understanding of the greater cycles of time and numbers, such as 1 habalatun = 460,800,000,000 years (Peat, 2002).

The legacy of this strategy is the still predominant belief that Western mathematics is the privileged manifestation of rationality of the human species and is the way that mathematics should look like (Goetzfridt, 2012). To exclude research of Indigenous mathematics from 'real mathematics', the terms pseudo-mathematics, proto-mathematics (Thomas, 1996) and ethnomathematics (Ascher, 2017; Ascher & Ascher, 1986) were coined.

The emergence of ethnomathematics has attracted controversies. Ethno-mathematicians believe that mathematics is culturally specific and each cultural group develops its own ways and styles of explaining, understanding and living with the help of their mathematical knowledge (Ascher, 2017). However, this focus on 'culture-specific mathematics' has deepened the ignorance from mainstream mathematicians on any study of Indigenous mathematics and believe that these studies add nothing to mainstream mathematics (Vithal & Skovsmose, 1997). What's worse, ethnomathematics is regarded by many mainstream mathematicians as an academic counterpoint to the universalism and internationalism of mathematics (Pais, 2013).

Our responses to the Gatekeeper critiques are based on our view that mathematics is part of the cultural heritage of all peoples, as intrinsic to our humanity as art and probably as ancient. We need to acknowledge that mathematics is a pan-cultural phenomenon and has many faces (Millroy, 1991). The differences between cultures pertain to the ways in which mathematics is expressed. Mainstream mathematics is top-down and decontextualized (de Almeida, 2019). In contrast, mathematics is dispersive in all aspects of Indigenous people's lives. This is why culture has been put at the core of ethnomathematics and probably research on Indigenous mathematics in the future.

If you ask a people who doesn't understand the algebra and anything to work on a formula, it is a bit unfair as she has no idea in both the concept, ideas and world view behind that for such kind of a decontextualized mathematical ideas so called. We are not creating an overwhelming notion that everything is mathematics. But we need to acknowledge that mathematics is a pan-cultural phenomenon and has many faces (Millroy, 1991). For example, almost everyone today uses a base-10 system for representing integers, but few people could imagine that diverse Indigenous groups use systems ranging from base-2 to base-68 (Dwyer & Minnegal, 2016; Owens & Lean, 2018).

Indigenous knowledge may be one of the keys to understanding how best to add new insights and advance mainstream and Indigenous mathematics together, but the mathematics research and teaching communities have to earn the respect of Indigenous knowledge holders. With two-way learning, Indigenous mathematical students, teachers and researchers will continue to thrive and the growing awareness, recognition and value of this knowledge will empower the communities to whom these stories belong.

Many studies have proven that learning mathematics with an indigenous mathematical approach can improve students' interests in studying mathematics and their performance (Fernández-Oliveras et al., 2021; Roza et al., 2020). We are also happy and optimistic for the future to note that Aboriginal and Torres Strait Islander mathematical elaborations have now been included in the F-10 Australian schools curriculum (https://www.australiancurriculum.edu.au).

**Statements and Declarations**

There are no conflicts of interests to declare. This research is not funded.